\documentclass{amsart}
\usepackage{amssymb,amsmath,amsfonts,amsthm}

\newcommand{\as}{\subseteq^*}
\newcommand{\CH}{The Continuum Hypothesis}
\newcommand{\ch}{the Continuum Hypothesis}
\newcommand{\abs}[1]{\left| #1 \right|}

\newcommand{\be}{\begin{enumerate}}
\newcommand{\ee}{\end{enumerate}}
\newcommand{\itm}{\item}
\newcommand{\cF}{\mathcal{F}}
\newcommand{\cU}{\mathcal{U}}

\newcommand{\pcc}{point-cofinite cover}
\newcommand{\un}{\bigcup}
\newcommand{\sub}{\subseteq}

\newcommand{\fb}{\mathfrak{b}}

\newtheorem{thm}{Theorem}[section]
\newcommand{\bthm}{\begin{thm}} \newcommand{\ethm}{\end{thm}}
\newtheorem{prop}[thm]{Proposition}
\newcommand{\bprp}{\begin{prop}} \newcommand{\eprp}{\end{prop}}
\newtheorem{fact}[thm]{Fact}
\newcommand{\bfct}{\begin{fact}} \newcommand{\efct}{\end{fact}}
\newtheorem{qtn}[thm]{Question}
\newcommand{\bqtn}{\begin{qtn}} \newcommand{\eqtn}{\end{qtn}}
\newtheorem{prob}[thm]{Problem}
\newcommand{\bprb}{\begin{prob}} \newcommand{\eprb}{\end{prob}}
\newtheorem{lem}[thm]{Lemma}
\newcommand{\blem}{\begin{lem}} \newcommand{\elem}{\end{lem}}
\newtheorem{claim}[thm]{Claim}
\newcommand{\bclm}{\begin{claim}} \newcommand{\eclm}{\end{claim}}
\newtheorem{cor}[thm]{Corollary}
\newcommand{\bcor}{\begin{cor}} \newcommand{\ecor}{\end{cor}}
\newtheorem{conj}[thm]{Conjecture}
\newcommand{\bcnj}{\begin{conj}} \newcommand{\ecnj}{\end{conj}}
\theoremstyle{definition}
\newtheorem{defn}[thm]{Definition}
\newcommand{\bdfn}{\begin{defn}} \newcommand{\edfn}{\end{defn}}
\newtheorem{spec}[thm]{Specializing}
\newcommand{\bspc}{\begin{spec}} \newcommand{\espc}{\end{spec}}
\theoremstyle{remark}
\newtheorem{rem}[thm]{Remark}
\newcommand{\brem}{\begin{rem}} \newcommand{\erem}{\end{rem}}
\newtheorem{cnv}[thm]{Convention}
\newcommand{\bcnv}{\begin{cnv}} \newcommand{\ecnv}{\end{cnv}}
\newtheorem{exam}[thm]{Example}
\newcommand{\bexm}{\begin{exam}} \newcommand{\eexm}{\end{exam}}
\newcommand{\bpf}{\begin{proof}} \newcommand{\epf}{\end{proof}}

\title[Productively Lindel\"of spaces]{On productively Lindel\"of spaces}

\author[F. D. Tall]{Franklin D. Tall}
\address[Tall]{Department of Mathematics, University of Toronto, Toronto, Ontario, M5S 2E4, Canada}
\email{f.tall@utoronto.ca}
\thanks{The research of the first named author is partially supported by grant A-7354 of the
Natural Sciences and Engineering Research Council of Canada.}

\author[B. Tsaban]{Boaz Tsaban}
\address[Tsaban]{Department of Mathematics, Bar-Ilan University, Ramat-Gan 52900, Israel}
\email{tsaban@math.biu.ac.il}



\begin{document}

\begin{abstract}
The class of spaces such that their product with every Lindel\"of space is
Lindel\"of is not well-understood.
We prove a number of new results concerning such productively Lindel\"of spaces with
some extra property, mainly assuming \ch.
\end{abstract}

\maketitle

\section{Applications of elementary submodels}

A quick introduction to the method of elementary submodels in our context is given in the appendix.

\bdfn
  A topological space is \emph{productively Lindel\"of} if its product with every Lindel\"of space is Lindel\"of.
  A space is \emph{powerfully Lindel\"of} if its $\omega$-th power is Lindel\"of.
\edfn

\bprb[Michael \cite{Przymusinski1980, Przymusinski1984}]
Are productively Lindel\"of spaces powerfully Lindel\"of?
\eprb

\blem[Alster \cite{Alster1988}]\label{lem1}
\CH{} implies productively Lindel\"of $T_3$ spaces of weight $\le \aleph_1$ are powerfully Lindel\"of.
\elem

Since Lindel\"of first countable $T_2$ spaces have cardinality (and hence weight) at most continuum,
we see that, assuming \ch,
productively Lindel\"of first countable $T_3$ spaces are powerfully Lindel\"of \cite{AurichiTall}.
This can be extended, as follows.

\bthm\label{thm2}
\CH{} implies that productively Lindel\"of sequential $T_3$ spaces are powerfully Lindel\"of.
\ethm
\bpf
Let $X$ be productively Lindel\"of sequential $T_3$, and $\cU$ be an open cover of $X^\omega$.
Without loss of generality, $\cU$ is composed of basic open subsets of $X^\omega$.
Let $M$ be a countably closed elementary submodel of $H_\theta$ of size $2^{\aleph_0}=\aleph_1$, for $\theta$ regular and sufficiently large,
such that $M$ contains $X, \cU$, and anything else needed.
For any space $Y \in M$, let $Y_M$ be the topology on
$Y \cap M$ generated by the sets $U \cap M$ where $U\in M$ is open in $Y$.

Since $M$ is countably closed, $X\cap M$ is a closed subset of $X$ and thus also productively Lindel\"of.
Since every open set in $X_M$ is open in $X\cap M$ with the relative topology, $X_M$ is a continuous
image of $X\cap M$, and therefore $X_M$ is productively Lindel\"of.\footnote{A general form of this argument appears in
Junqueira-Tall \cite{Junqueira1998}, see Proposition \ref{p} in the appendix.}
The weight of $X_M$ is $\le |M|\le \aleph_1$, so by Alster's Lemma \ref{lem1}, $(X_M)^\omega$ is Lindel\"of.

Since $M$ is countably closed, $X^\omega\cap M = (X\cap M)^\omega$, that is, $(X^\omega)_M = (X_M)^\omega$ as sets.
Also, for $B_i \in M$ open in $X$, the set $\prod_{i<\omega} B_i\sub X^\omega$, when intersected with $M$, is just
$\prod_{i < \omega} (B_i \cap  M)$, so we see that as spaces, $(X^\omega)_M = (X_M)^\omega$.

Thus, $(X^\omega)_M$ is Lindel\"of.
As $\cU\in M$, we have by elementarity that $\{U \cap M : U \in \cU \cap M\}$ is an open cover of $(X^\omega)_M$.
Thus, there are $\{ U_n : n < \omega \} \sub \cU \cap M$ such that $\{ U_n \cap M : n < \omega \}$ covers $(X^\omega)_M$.
Since $M$ is countably closed, $\{ U_n : n < \omega \} \in M$.
$M \models \{U_n : n < \omega \}$ covers $X ^ \omega$, so indeed $\{U_n : n < \omega\}$ covers $X ^ \omega$.
\epf

It would be nice to eliminate the hypothesis that $X$ be sequential.
This was only used to get that the sequentially closed set $X\cap M$ is indeed closed.
There are large compact (hence productively Lindel\"of) sequential $T_3$ spaces,
so ``sequential'' is indeed an improvement over ``first countable''.
We do not know whether \ch{} is necessary for these results.

It is not known whether the weight restriction in Alster's Lemma \ref{lem1} can be removed, nor whether \ch{} is necessary.
There is no reason to believe the weight of a productively Lindel\"{o}f space cannot exceed its cardinality, so the following result is not obvious.

\bcor
\CH{} implies productively Lindel\"{o}f $T_3$ spaces of cardinality $\aleph_1$ are powerfully Lindel\"{o}f.
\ecor
\bpf
As in the proof of Theorem \ref{thm2}, we take a countably closed elementary submodel $M$ of $H_\theta$ with $|M|=\aleph_1$, such
that $M$ contains everything needed, and, \emph{in addition}, $X\sub M$.
It follows that the weight of $X_M$ is $\leq \aleph_1$.  Since $X \cap M = X$, $X_M$ is a continuous image of $X$ and hence is productively Lindel\"{o}f and so powerfully Lindel\"{o}f, and hence, as before, $X$ is powerfully Lindel\"{o}f.
\epf

\section{Selective covering properties}

\bdfn
A \emph{\pcc} of a space is an infinite open cover such that each point is in all but finitely many members of the cover.
\edfn

\bdfn\label{AHM}
A topological space $X$ is:
\be
\itm \emph{Alster} if each cover of $X$ by $G_\delta$ sets such that each compact set is included in one of them has a countable subcover.\footnote{Alster's terminology \cite{Alster1988} is slightly different, but equivalent.}
\itm \emph{Hurewicz} if for each sequence $\{\cU_n\}_{n<\omega}$ of open covers without finite subcovers, there are
finite $\cF_n\sub \cU_n$ such that $\{\un \cF_n : n < \omega\}$ is a \pcc.
\itm \emph{Menger} if for each sequence $\{\cU_n\}_{n<\omega}$ of open covers without finite subcovers, there are
finite $\cF_n\sub \cU_n$ such that $\{\un \cF_n : n < \omega\}$ is a cover.
\ee
Let P be a property of topological spaces. A space $X$ is \emph{powerfully P} if $X^\omega$ has the property P.
$X$ is \emph{finitely powerfully P} if all its finite powers of $X$ have the property P.
\edfn

By the definition, Hurewicz spaces are Menger.
Alster spaces are productively and powerfully Lindel\"of \cite{Alster1988}.
A slightly extended version of an argument from \cite{AurichiTall} yields the following generalization of results
from \cite{AurichiTall, Tall}.

\bthm\label{AH}
Alster spaces are Hurewicz.
\ethm
\bpf
Let $\{\cU_n\}_{n<\omega}$ be a sequence of open covers of $X$ without finite subcovers.
We may assume that each $\cU_n$ is closed under finite unions.
Let
$$\cU=\left\{\bigcap_{n<\omega} U_n : \forall n,\ U_n\in\cU_n\right\}.$$
Since $X$ is Alster, there is a countable subcover $\{V_m : m<\omega\}$ of $\cU$.
For each $m$, write $V_m=\bigcap_{n<\omega} U_{mn}$, where $U_{mn}\in\cU_n$ for all $n$.
Then for each $x\in X$, $x\in \un_{m\le n}U_{mn}$ for all but finitely many $n$.
\epf

%\brem
%The proof of Theorem \ref{AH} shows that, if $X$ is Alster, then for each sequence
%$\{\cU_n\}_{n<\omega}$ of covers by $G_\delta$ sets such that for each $n$, $\cU_n$ has no finite subcover of $X$,
%but for each compact subset of $X$ there is a finite subcover in $\cU_n$, then
%there are finite $\cF_n\sub \cU_n$ such that $\{\un \cF_n : n < \omega\}$ is a \pcc{} of $X$.
%\erem

Thus, each property in Definition \ref{AHM} implies the next one. These implications
are strict: For sets of reals (indeed, for arbitrary spaces where every compact set is $G_\delta$),
Alster is clearly equivalent to \emph{$\sigma$-compact},
and Hurewicz fits strictly between $\sigma$-compact and Menger.\footnote{For an accessible exposition of this result, see \cite{1T}.}

\bcor
Alster spaces are finitely powerfully Hurewicz.
\ecor
\bpf
Finite products of Alster spaces are Alster \cite{Alster1988, Barr2007}.
The proof in \cite{Barr2007} does not use separation axioms.
Apply Theorem \ref{AH}.
\epf

A set of reals is \emph{totally imperfect} if it includes no uncountable perfect (equivalently, compact) set.
\bthm
There is a finitely powerfully Hurewicz set of reals which is not productively Lindel\"{o}f (and hence not Alster).
\ethm
\bpf
Michael \cite{Michael1971} proved that totally imperfect set of reals are not productively Lindel\"{o}f.
Bartoszy\'nski and the second named author \cite{1BT} proved that there is a totally imperfect, finitely powerfully Hurewicz set of reals.
\epf

\blem[Alster \cite{Alster1988}]\label{lem6}
  \CH{} implies productively Lindel\"of $T_3$ spaces of weight $\le \aleph_1$ are Alster.
\elem

Alster asked whether every productively Lindel\"of space is Alster \cite{Alster1988}.
Alster's problem is still open. The following problem may be easier.

\bprb
  Does \ch{} imply productively Lindel\"of sequential $T_3$ spaces are Alster?
\eprb

\bthm\label{newThm14}
  \CH{} implies productively Lindel\"of sequential $T_3$ spaces are finitely powerfully Hurewicz.
\ethm

\bpf
  It suffices to show that if $(X^k)_M$ is Hurewicz, then $X^k$ is Hurewicz.
  We assume without loss of generality that the sequence
  $\{\cU_n\}_{n < \omega}$ of open covers is in $M$.
  Then for each $n$, there is a finite $\cF_n \sub \cU_n$ such that
  $\bigcup \{ U \cap M : U \in \cF_n \}$ is a \pcc{} of $(X^k)_M$.
  Note that each $\cU_n \in M$. Since $\cU_n$ is countable, it is included in $M$, and hence each $\cF_n \in M$.
  As $M$ is countably closed, $\{ \cF_n \}_{n < \omega} \in M$. Thus,
$$M \vDash \left\{\un\cF_n : n<\omega\right\}\mbox{ is a \pcc{} of }X^k,$$
and therefore the same holds in ``the real world'', so indeed $X^k$ is Hurewicz.
\epf

D-spaces were defined in \cite{vanDouwenPfeffer}. See also \cite{Eisworth2007, Gruenhage2009}.

\bdfn
  A space $X$ is \emph{D} if for every \textit{neighborhood assignment}
  $\{V_x\}_{x\in X}$ (i.e., each $V_x$ is an open set containing $x$),
  there is a closed discrete $Y\sub X$ such that $\{V_x\}_{x\in
    Y}$ covers $X$.
\edfn

Aurichi \cite{Aurichi} proved that Menger spaces are D. Thus, assuming \ch{},
productively Lindel\"of sequential $T_3$ spaces are finitely productively D.
L. Zdomskyy pointed out to us that this last assertion can be generalized substantially.
A \emph{Michael space} is a Lindel\"of space $M$ such that $M \times \mathbb{P}$ (the space of irrationals) is not Lindel\"of.
Michael spaces can be constructed from a variety of axioms (in particular, from \CH{}),
and it is a major open problem whether they can be constructed outright in ZFC.
If there is a Michael space $M$, then productively Lindel\"of spaces are Menger (and thus D) \cite{1RZ}.
Indeed, Zdomskky proves in \cite{SF1} that if $X$ is not Menger, then $\mathbb{P}$ is a compact-valued upper-semicontinuous image of $X$.
Thus, if $X$ is not Menger, then the non-Lindel\"of space $M\times\mathbb{P}$ is a compact-valued upper-semicontinuous image of $M\times X$. 
Consequently, $M\times X$ is not Lindel\"of.

\section{Indestructibly productively Lindel\"of spaces}

\bdfn
A space is \emph{indestructibly productively Lindel\"of} if it is productively
Lindel\"of in every countably closed forcing extension.
\edfn

Aurichi and the first named author proved that a metrizable space is indestructibly
productively Lindel\"of if and only if it is $\sigma$-compact \cite{AurichiTall}.
It is easily seen that if a space $Y$ is Hurewicz in a countably closed
extension, then it is Hurewicz.
The following theorem answers a question of Aurichi and the first named author \cite{AurichiTall}.

\bthm\label{thm7}
Indestructibly productively Lindel\"of $T_3$ spaces are powerfully Lindel\"of and finitely powerfully Hurewicz (in particular, finitely powerfully D).
\ethm
\bpf
Powerfully Lindel\"of: Collapse $\max(2^{\aleph_0}, w(X))$ to $\aleph_1$ via countably
closed forcing. In the extension, the indestructibly productively Lindel\"of $X$ remains
productively Lindel\"of and hence, by \ch, $X$ becomes Alster. Then by
Lemma \ref{lem1}, $X$ becomes powerfully Lindel\"of. But as a set, $X^\omega$
in the extension is the  same as $X^\omega$ in the ground model. Since
the space $X^\omega$ is Lindel\"of in a countably closed extension, it is
Lindel\"of in the ground model, as claimed.

Finitely powerfully Hurewicz: In the extension obtained by collapsing as above,
$X$ is Alster, so every finite power of $X$ is Hurewicz.
But then every finite power of $X$ is Hurewicz in the ground model.
\epf

\bprb
Are indestructibly productively Lindel\"{o}f spaces Alster?
\eprb

\bcor
Indestructibly productively Lindel\"of $p$-spaces are $\sigma$-compact.
\ecor
\bpf
Let $X$ be indestructibly productively Lindel\"of and a $p$-space in the sense
of Arhangel'ski\u\i{} \cite{Arhangelskii1963}. Then, as a
paracompact $p$-space, $X$ maps perfectly onto a metrizable $Y$. Let $Z$ be
Lindel\"of in a countably closed extension. Then $X \times Z$ is Lindel\"of. But then, by continuity, so is
$Y \times Z$. So $Y$ is indestructibly productively Lindel\"of. But for metrizable
spaces, indestructible productive Lindel\"ofness is equivalent to
$\sigma$-compactness \cite{AurichiTall}, which latter property is a perfect
invariant.
\epf

Arhangel'ski\u{\i} \cite{Arhangelskii1986} proved that if $X^\omega$ is Lindel\"of,
then either $X$ is compact or $X^\omega$ includes a closed copy of $\mathbb{P}$.
Since the latter option is impossible for Menger spaces, he concluded that
$X$ is powerfully Menger if and only if $X$ is compact.\footnote{Arhangel'ski\u\i{} denotes by ``Hurewicz" the property we call ``Menger".
We use the currently accepted terminology.}

\bcor
If there is a Michael space, then for every space $X$,
$X^\omega$ is productively Lindel\"of if and only if $X$ is compact.\qed
\ecor

\bthm
If $X^\omega$ is indestructibly productively Lindel\"of, then $X$ is compact.
\ethm
\bpf
Again, collapse $\max(2^{\aleph_0}, w(X))$ to $\aleph_1$. In the extension,
$X^\omega$ is productively Lindel\"of and there is a Michael space, since \ch{} holds, which implies that there is a Michael space \cite{Michael1971}. Therefore $X$ is compact in the extension, and so is compact.
\epf

\section{Mengerizing Michael's problems}

As Menger implies Lindel\"of, the classic problems about productively
Lindel\"of spaces make sense when \emph{Lindel\"of} is replaced by \emph{Menger}.

\bexm
$\mathbf{\omega}$ (the countable discrete space) is productively Menger, but not powerfully Menger.
\eexm

The product of a Menger space with $\mathbb{P}$ cannot be Menger since $\mathbb{P}$ is not Menger,
but the question whether the product of a Menger space with $\mathbb{P}$ must be Lindel\"{o}f is less trivial.
We will show that the answer is negative, in a very strong sense.

\bdfn
An open cover $\cU$ of a space $X$ is an \emph{$\omega$-cover} if $X\notin\cU$, but for each finite
subset $F$ of $X$ there is $U\in\cU$ containing $F$. $X$ is a \emph{$\gamma$-space} if each $\omega$-cover
of $X$ includes a \pcc.
\edfn

$\gamma$-spaces were introduced by Gerlits and Nagy \cite{1GN}, who proved that, for Tychonoff spaces,
$X$ is a $\gamma$-space if and only if the space $C_p(X)$ (the continuous the real-valued functions on $X$
with the topology of pointwise convergence) is Fr\'echet-Urysohn. This is a very strong property.
It is, for example, consistent that all metrizable $\gamma$-spaces are countable \cite{1GN}. If $X$ is a $\gamma$-space
then $X$ is Hurewicz. Being a $\gamma$-space is preserved by finite powers \cite{1GN}. In particular,
$\gamma$-spaces are finitely powerfully Hurewicz.

For $f,g\in\omega^\omega$, $f\le^* g$ means that $f(n)\le g(n)$ for all but finitely many $n$.
A subset of $\omega^\omega$ is \emph{unbounded} if it is unbounded with respect to $\le^*$.
The minimal cardinality of an unbounded subset of $\omega^\omega$ is denoted $\fb$.
$\aleph_1\le\fb\le 2^{\aleph_0}$. In particular, \ch{} implies $\fb=\aleph_1$.
Additional information on $\fb$ and similar combinatorial cardinal characteristics of the continuum can be found
in \cite{vanDouwen}.

We identify elements $x\in [\omega]^\omega$ with increasing elements of $\omega^\omega$
by letting $x(n)$ be the $n$th element of $x$.
We will need the following well-known fact. For the reader's convenience, we reproduce here
the proof given in \cite{1TO}.

\blem[folklore]\label{infdisj}
If $B\sub[\omega]^\omega$ is unbounded, then for each increasing $f\in\omega^\omega$,
there is $x\in B$ such that $x\cap (f(n),f(n+1))=\emptyset$ for infinitely
many $n$.
\elem
\bpf
Assume that $f$ is a counterexample. Let $g$ dominate all functions $f_m(n)=f(n+m)$, $m<\omega$.
Then for each $x\in B$, $x\le^* g$. Indeed, let $m$ be such that for all $n\ge m$,
$x\cap (f(n),f(n+1))\neq\emptyset$. Then for each $n$, the $n$-th element of $x$
is smaller than $f_{m+1}(n)$.
\epf

Orenshtein and the second named author \cite{1TO} proved that an assumption weaker than $\fb=\aleph_1$ implies that
there is an uncountable $\gamma$-space $X\sub\mathbb{R}$.
The proof of the forthcoming Theorem \ref{main} is a modification of their proof, slightly simplified in light of the stronger assumption.

\bdfn
Identify $P(\omega)$ with the Cantor space $2^\omega$, using characteristic functions.
This defines the \emph{Cantor topology} on $P(\omega)$.
Consider the finer, \emph{Michael topology} on $P(\omega)$ obtained by declaring all elements of $[\omega]^\omega$
isolated. Henceforth, unless otherwise indicated, $P(\omega)$ is always considered with the Michael topology.
\edfn

The basic open sets in the Cantor topology of $P(\omega)$ are thus those of the form
$$[s,n]=\{x\in P(\omega) : x\cap\{0,\dots,n-1\}=s\},$$
where $n\in\omega$ and $s\sub \{0,\dots,n-1\}$.
We will use the following modification of Lemma 1.2 of Galvin and Miller \cite{GM84}.

\blem\label{GM+}
Consider $P(\omega)$ with the Michael topology.
Assume that $[\omega]^{<\omega}\sub Y\sub P(\omega)$, $Y$ is countable, and $\cU$ is a family of open subsets of $P(\omega)$ such
that each finite subset of $Y$ is included in some member of $\cU$.
There are $m_0<m_1<\dots$ and (not necessarily distinct) $U_0,U_1,\dots\in\cU$ such that:
\be
\itm For each $y\in Y$, $y\in U_n$ for all but finitely many $n$.
\itm For each $x\sub\omega$, $x\in U_n$ whenever $x\cap (m_n,m_{n+1})=\emptyset$.
\ee
\elem
\bpf
Enumerate $Y=\{y_n : n<\omega\}$.

Let $m_0=0$. For each $n\ge 0$: Take $U_n\in\cU$, such that $P(\{0,\dots,\allowbreak m_n\})\cup\{y_0,\dots,y_n\}\sub U_n$.
Let $s\sub \{0,\dots,m_n\}$. As $s\in [\omega]^{<\omega}$ and $U_n$ is a neighborhood of $s$,
$U_n$ includes a neighborhood of $s$ in the Cantor set topology, and thus there is $k_s$ such that
for each $x\in P(\omega)$ with $x\cap\{0,\dots,k_s-1\}=s$, $x\in U_n$. Let $m_{n+1}=\max\{k_s : s\sub\{0,\dots,m_n\}\}$.
\epf

Since $P(\omega)$ is equipped with a topology finer than that of Cantor's space, which is metrizable,
the following result cannot be proved outright in ZFC.

\bthm\label{main}
Consider $P(\omega)$ with the Michael topology.
If $\fb=\aleph_1$, then there is a $\gamma$-space $X\sub P(\omega)$ which is also a Michael space (i.e., such that $X\times\mathbb{P}$ is not Lindel\"of).
\ethm
\bpf
For $x,y\in [\omega]^\omega$, $x\as y$ means that $x\setminus y$ is finite.
As $\fb=\aleph_1$, there is an unbounded (with respect to $\le^*$)
set $\{x_\alpha : \alpha<\aleph_1\}\sub[\omega]^\omega$
such that for all $\alpha<\beta<\aleph_1$, $x_\beta\as x_\alpha$.

Let $$X=\{x_\alpha : \alpha<\aleph_1\}\cup[\omega]^{<\omega}\sub P(\omega),$$
with the subspace topology (so that the elements $x_\alpha$ are isolated),
and consider $X\times [\omega]^\omega$, where the space $[\omega]^\omega$ on the right is endowed with the
ordinary Cantor space topology, so that it is homeomorphic to $\mathbb{P}$.
The uncountable set $\{(x_\alpha,x_\alpha) : \alpha<\aleph_1\}$ is closed and discrete in $X\times [\omega]^\omega$.
Thus, this space is not Lindel\"of.
Once we prove that $X$ is a $\gamma$-space, we will have in particular that $X$ is Lindel\"of, so that $X$ is a Michael space.
That $X$ is a Michael space is essentially proved in \cite{vanDouwen}; that $X$ is a $\gamma$-space is new.

Let $\cU$ be an $\omega$-cover of $X$.
For each $\alpha<\aleph_1$, let $X_\alpha=\{x_\beta : \beta<\alpha\}\cup[\omega]^{<\omega}$.

Let $\alpha_0=0$. By Lemma \ref{GM+}, there are $m^0_0<m^0_1<\dots$ and elements $U^0_0,U^0_1,\dots\in\cU$
such that each member of $X_{\alpha_0}$ is in $U^0_n$ for all but finitely many $n$, and for each
$x\in P(\omega)$, $x\in U^0_n$ whenever $x\cap (m^0_n,m^0_{n+1})=\emptyset$.
Let $D_0=\omega$.

As $\alpha_0<\aleph_1$, $\{x_\alpha : \alpha_0<\alpha<\aleph_1\}$ is unbounded.
By Lemma \ref{infdisj}, there is $\alpha_1>\alpha_0$ such that
$D_1=\{n : x_{\alpha_1}\cap (m^0_n,m^0_{n+1})=\emptyset\}$ is infinite.
By Lemma \ref{GM+}, there are $m^1_0<m^1_1<\dots$ and members $U^1_0,U^1_1,\dots\in\cU$
such that each member of $X_{\alpha_1}$ is in $U^1_n$ for all but finitely many $n$, and for each
$x\in P(\omega)$, $x\in U^1_n$ whenever $x\cap (m^1_n,m^1_{n+1})=\emptyset$.

Continue in the same manner to define, for each $k>0$, elements with the following properties:
\be
\itm $\alpha_k>\alpha_{k-1}$;
\itm $D_k=\{n : x_{\alpha_k}\cap (m^{k-1}_n,m^{k-1}_{n+1})=\emptyset\}$ is infinite;
\itm $m^k_0<m^k_1<\dots$;
\itm $U^k_0,U^k_1,\dots\in\cU$;
\itm each member of $X_{\alpha_k}$ is in $U^k_n$ for all but finitely many $n$; and
\itm For each $x\in P(\omega)$, $x\in U^k_n$ whenever $x\cap (m^k_n,m^k_{n+1})=\emptyset$.
\ee
Let $\alpha=\sup_k\alpha_k$. Then $\alpha<\aleph_1$, $X_\alpha$ is countable, and $X_{\alpha_k}\sub X_{\alpha_{k+1}}$ for all $k$.
Thus, there are for each $k$ a finite $F_k\sub X_{\alpha_k}$
such that $F_k\sub F_{k+1}$ for all $k$, and $X_\alpha=\un_kF_k$.
For each $k$, let $I_k = \{n\in D_k : F_k\sub U^k_n\}$.
$I_k$ is an infinite (indeed, cofinite) subset of $D_k$, and
for each $x\in X_\alpha$, if $N$ is the first with $x\in F_N$,
then $x$ belongs to  $\bigcap_{n\in I_k}U^k_n$ for all $k\ge N$.

Take $n_0\in I_1$.
For $k>0$, take $n_k\in I_{k+1}$ such that $m^k_{n_k}>m^{k-1}_{n_{k-1}+1}$,
$x_\alpha\cap (m^k_{n_k},m^k_{n_k+1})\sub x_{\alpha_{k+1}}\cap (m^k_{n_k},m^k_{n_k+1})$,
and $U^k_{n_k}\notin\{U^1_{n_1},\dots,U^{k-1}_{n_{k-1}}\}$.
We claim that each member of $X$ is in $U^k_{n_k}$ for all but finitely many $k$.\footnote{Technically,
\pcc{}s are required to be infinite. To see that this follows, note that if $\{U^k_{n_k} : k<\omega\}$ is finite, then
there is $U$ such that $U=U^k_{n_k}$ for infinitely many $k$. As $X$ is not in $\cU$,
there is $x$ which is not in $U$, and consequently not in infinitely many members of the
sequence $\{U^k_{n_k}\}_{k<\omega}$, contradicting the assertion we are about to prove.}
By the last assertion in the previous paragraph, this is true for each member of
$X_\alpha$. As for each $\beta\ge\alpha$ we have that $x_\beta\as x_\alpha$, it suffices to
show that for each
$x\as x_\alpha$, $x\in U^k_{n_k}$ for all but finitely many $k$.
For each large enough $k$, $m^{k}_{n_k}$ is large enough, so that
$$x\cap (m^k_{n_k},m^k_{n_k+1})\sub x_\alpha\cap (m^k_{n_k},m^k_{n_k+1})\sub
x_{\alpha_{k+1}}\cap (m^k_{n_k},m^k_{n_k+1})=\emptyset,$$
since $n_k\in D_{k+1}$. Thus, $x\in U^k_{n_k}$.
\epf

\section{Analytic spaces}

The first named author proved in \cite{Tall} that
every analytic, metrizable, productively Lindel\"of space is $\sigma$-compact
if and only if there is a Michael space.
The hypothesis of metrizability can be removed.
According to Arhangel'ski\u{\i} \cite{Arhangelskii1986},
a space is \emph{analytic} if it is a continuous image of the space
$\mathbb{P}$ of irrationals.

\bthm
Every analytic productively Lindel\"of space is $\sigma$-compact if and only
if there is a Michael space.
\ethm
\bpf
By their definition, analytic spaces are Lindel\"of. Relying on results of Jayne
and Rogers \cite{JayneRogers}, Arhangel'ski\u{\i} \cite{Arhangelskii1986} proved that
analytic spaces are perfect pre-images of metrizable spaces.
Since both productive Lindel\"ofness and $\sigma$-compactness are perfect
invariants, we are done.
\epf

Perfect pre-images of analytic spaces are called \emph{properly analytic} in
\cite{JayneRogers}. It follows immediately that every properly analytic,
productively Lindel\"of space is $\sigma$-compact if and
only if there is a Michael space.

According to Hansell \cite{Hansell}, a space is \emph{$K$-analytic} if
it is the continuous image of a Lindel\"of \v{C}ech-complete space.

\bprb
Is it consistent that every productively Lindel\"of $K$-analytic space is
$\sigma$-compact?
\eprb

The first named author also proved in \cite{Tall} that
the Axiom of Projective Determinacy implies that every projective,
metrizable, productively Lindel\"of space is $\sigma$-compact if and only if
there is a Michael space.

We can certainly extend this to perfect pre-images of
projective metrizable spaces, but what is the analog of
Arhangel'ski\u{\i}'s definition? One possibility is to define
``projective'' as a continuous image of a projective subset of
$\mathbb{P}$ (or $\mathbb{R}$). We do not know whether this definition
allows us to apply Projective Determinacy as desired. However, we do
have the following.

\bthm\label{Thm14}
\CH{} implies every productively Lindel\"of, continuous
image of a separable metrizable space is $\sigma$-compact.
\ethm
\bpf
Let $X$ be such an image. $X$ is $T_3$. $X$ has a countable network and so $X$
is separable and every closed subset is $G_\delta$.
The weight of $X$ is $\le 2^{\aleph_0}$, so by Alster's Lemma \ref{lem6},
$X$ is Alster. But Alster spaces in which compact sets are $G_\delta$ are
$\sigma$-compact \cite{Alster1988}.
\epf

The \emph{Baire Hierarchy} is formed by closing the collection of
\emph{closed} sets under countable unions and intersections. In contrast
to the Borel Hierarchy, the Hurewicz Dichotomy fails at a low level. A
\emph{$K_{\sigma \delta}$}\textit{ space} is a space which is the intersection of
countably many $\sigma$-compact subspaces of some larger space.

\bexm\label{ex}
There is a $K_{\sigma \delta}$ space which is neither
$\sigma$-compact nor includes a closed copy of $\mathbb{P}$.
\eexm
\bpf
In \cite{Arhangelskii2000} Arhangel'ski\u\i{} constructs a $K_{\sigma \delta}$
space, due to Okunev, which is not $\sigma$-compact but has only one
non-isolated point, so does not include a closed copy of $\mathbb{P}$.
The space is obtained by taking the Alexandrov duplicate of $\mathbb{P}$,
and then collapsing the non-discrete copy of $\mathbb{P}$ to a point.
\epf

\subsection{$\kappa$-analytic spaces}\label{lub}
Descriptive set-theorist Ben Miller told us that the ``right''
definition of projective in a non-separable metrizable context is
the following one.

\bdfn
  A $T_2$ space $X$ is \emph{$\kappa$-analytic}, where $\kappa$ is an
  uncountable cardinal, if $X$ is a continuous image of the product of
  $\aleph_0$ copies of the discrete space of size $\kappa$.
\edfn

Every space $X$ is $\abs{X}$-analytic.

Recall that, according to the Hurewicz Dichotomy,
every analytic non-$\sigma$-compact subspace of the
Baire space contains a closed copy of the Baire space, and thus,
if there is a Michael space, an analytic metrizable space is productively Lindel\"of
if and only if it is $\sigma$-compact.
This and Example \ref{ex} motivate the following question.

\bqtn\label{qn}
Let $\kappa<2^{\aleph_0}$. Is it consistent that
\be
\itm Every non-$\sigma$-compact $\kappa$-analytic metrizable space includes a closed copy of $\mathbb{P}$?
\itm Every productively Lindel\"of $\kappa$-analytic metrizable space is $\sigma$-compact?
\ee
\eqtn

It was L. Zdomskyy who pointed out to us that, if we drop the metrizability assumption,
then the one-point Lindel\"{o}fication of the discrete space of size $\aleph_1$
gives a $T_3$ counter-example to both items of Question \ref{qn}, and that
$\kappa < 2^{\aleph_0}$ is necessary for the problem to have a possibly genuine descriptive set theoretic flavor.

The hypotheses in the following theorem, which answers (2) of Question \ref{qn},
follow from Martin's Axiom plus the negation of \ch{}, see \cite{1A}.

\bthm[Zdomskyy]
Assume that there is no cover of the Cantor space by $\aleph_1$ meager sets, and there is a Michael space.
Then every productively Lindel\"{o}f $\aleph_1$-analytic subset of the Cantor space is $\sigma$-compact.
\ethm
\bpf
The proof is analogous to that of Theorem 3.3 in \cite{1RZ}.
The key is to observe that the proof of Repicky's Theorem in \cite{1R}, that
$\mathbf{\Sigma}_2^1$ subspaces of $\mathbb{P}$ that are not the union of
$\aleph_1$ compact subspaces include a closed copy of $\mathbb{P}$, also works for
$\aleph_1$-analytic subspaces.
\epf

\bprb[Zdomskyy]
Is it consistent that every non-$\sigma$-compact space which is
$\kappa$-analytic for some uncountable cardinal $\kappa < 2^{\aleph_0}$ includes a closed copy of $\mathbb{P}$?
\eprb

\section{Spaces of countable type}

In 1957, M. Henriksen and J. Isbell \cite{HenriksenIsbell}
introduced the class of
(Tychonoff) spaces that are \emph{Lindel\"of at infinity}, i.e., the
complement $\beta X\setminus X$ of the space $X$ in its Stone-\v{C}ech
compactification is Lindel\"of. They proved that
a Tychonoff space $X$ is Lindel\"of at infinity if and only if each
compact subset of $X$ is included in a compact $K \sub X$ such that
$\chi(X,K) \le \aleph_0$, i.e., there is a countable base for the neighborhoods
of $K$ in $X$.
Arhangel'ski\u\i{} \cite{Arhangelskii1965} called spaces satisfying the latter
equivalent condition \emph{of countable type}. Locally compact spaces,
metrizable spaces, \v{C}ech-complete spaces, and their common generalization,
$p$-spaces, are all of countable type.

We present a simple proof for the following generalization of
a result of Alster from \cite{Alster1988}.

\bthm[Alas, et al.\ \cite{AAJT}]\label{thm29}
\CH{} implies every productively Lindel\"of $T_3$ space
of countable type and weight $\le \aleph_1$ is $\sigma$-compact.
\ethm
\bpf
We generalize Michael's original proof that \ch{} implies productively Lindel\"of metrizable spaces are
$\sigma$-compact.

Embed $X$ in $[0,1]^{\aleph_1}$.
Its closure in $[0,1]^{\aleph_1}$ is a compactification $\gamma X$ of $X$.
The identity map on $X$ extends to a continuous surjection $f:\beta X\to\gamma X$,
and since it fixes $X$, $f$ maps $\beta X \setminus X$ onto $\gamma X \setminus X$
\cite[3.5.7]{Engelking1989}.
As $\beta X \setminus X$ is Lindel\"of, so is $\gamma X \setminus X$.

Assume that $X$ is not $\sigma$-compact. Then $\gamma X\setminus X$ is not $G_\delta$ in $\gamma X$. 
By \ch, we can take a collection $\{ U_\alpha : \alpha < \aleph_1\}$ of open sets including
$\gamma X \setminus X$, such that every open set including $\gamma X \setminus X$ includes some
$U_\alpha$. By taking countable intersections and thinning out, we can
find a strictly decreasing sequence $\{ V_\beta \}_{\beta < \aleph_1}$
of $G_\delta$ sets including $\gamma X \setminus X$, such that every open set including
$\gamma X \setminus X$ includes some $V_\beta$. For each $\beta < \aleph_1$, take $p_\beta \in (V_{\beta + 1} \setminus V_\beta)\cap X$. 

Let 
$$Y=(\gamma X \setminus X) \cup \{ p_\beta : \beta <\aleph_1 \}$$
Put a topology on $Y$ by strengthening the subspace topology to make each $\{ p_\beta \}$ open. 
The usual Michael space argument \cite{Michael1971} (cf.\ \cite{AAJT} for more general arguments)
shows that $Y$ is Lindel\"of, but its product with $X$ is not.
Indeed, 
$Y$ is Lindel\"of, since each open set including $\gamma X\setminus X$ includes all but countably many $p_\beta$'s.
To see that $Y \times X$ is not Lindel\"of, note that the set $\{(p_\beta, p_\beta) : \beta < \aleph_1\}$
is closed and discrete in $Y \times X$.
\epf

$X$ is \emph{absolute Borel} if it is Borel in $\beta X$.
In this case, $\beta X \setminus X$ is a Baire subspace of $\beta X$.
As Baire subspaces of compact $T_2$ spaces are Lindel\"of \cite{HIJ},
Lindel\"of absolute Borel spaces are of countable type. We therefore have the following.

\bcor
\CH{} implies productively Lindel\"of absolute Borel spaces
of weight $\le \aleph_1$ are $\sigma$-compact.
\ecor

\subsection*{Acknowledgments}
We thank Ofelia Alas for correcting the original argument of the proof of Theorem \ref{thm29},
and the referee for useful comments and suggestions.
We owe special thanks to Lyubomyr Zdomskyy for inspiring discussions and suggestions leading to
some of the results in this paper. In particular, we thank him for his contributions in Subsection \ref{lub},
which clarify a question from an earlier version.

\appendix

\section{Some remarks on elementary submodels and forcing}
For the reader not so familiar with elementary submodels, we make some elementary remarks which may be helpful in understanding the proofs
in this paper which involve this method.

First of all, the sets $H_\theta$ in our proofs
appear only for technical reasons; we really think instead of the universe $V$.  For elucidation of this point, see Chapter 24 of \cite{1JW}.

An elementary submodel $M$ is \textit{countably closed} if each countable subset of $M$ is a member of $M$.
For such models $M$, if $X \in M$, then the collection of countable sequences of members of $X \cap M$ is
the same as the collection of countable sequences lying in $M$ of members of $X$, i.e., $(X\cap M)^\omega=X^\omega\cap M$.
A straightforward closing-off (L\"{o}wenheim-Skolem) argument establishes that $H_\theta$, for regular $\theta \geq 2^{\aleph_0}$, has a countably closed elementary submodel of size $2^{\aleph_0}$.

\bdfn
For a topological space $X$ with topology $\tau$, $X_M$ is the topological space $X\cap M$ with the topology with basis
$\{U\cap M : U\in\tau\cap M\}$.
\edfn
The proofs of the following basic facts are illustrative.

\blem[folklore]\label{sc}
Assume that $X$ is a $T_2$ space, $M$ is a countably closed elementary submodel of $H_\theta$ for some sufficiently large regular $\theta$, and
$X\in M$. Then: $X\cap M$ is a sequentially closed subset of $X$.
\elem
\bpf
Let $\{x_n\}_{n<\omega}$ be a sequence of elements of $X\cap M$
converging to a point $x\in X$. As $M$ is countably closed and each $x_n\in M$, $\{x_n\}_{n<\omega}\in M$ as well.
By elementarity, $M\vDash \{x_n\}_{n<\omega}$ converges to some point $y$. Since $X$ is Hausdorff, $x=y\in M$, so $x\in X\cap M$.
\epf

\bprp[Junqueira-Tall \cite{Junqueira1998}]\label{p}
Assume that:
\be
\itm[(a)] $X$ is a sequential $T_2$ space;
\itm[(b)] $M$ is a countably closed elementary submodel of $H_\theta$ for some sufficiently large regular $\theta$; and
\itm[(c)] $X\in M$.
\ee
Then:
\be
\itm $X\cap M$ is a closed subset of $X$.
\itm $X_M$ is a continuous image of $X\cap M$.
\itm For each property $P$ of $X$ preserved by continuous images and closed subspaces, $X_M$ has the property $P$.
\ee
\eprp
\bpf
(1) follows from Lemma \ref{sc}, as $X$ is sequential.

(2) The identity map from $X\cap M$ with the relative topology onto $X_M$ is continuous, since every open set in $X_M$ is
open in $X\cap M$.

(3) follows from (1) and (2).
\epf

Another observation about countably closed models is that, roughly speaking, if properties involving countable sets
(such as Lindel\"{o}fness) are true for the fragment of $X$ lying in $M$, then $M$ will demonstrate that,
and thus, by elementarity, $X$ will really have that property.
Thus, such properties as powerfully Lindel\"{o}f, (finitely) powerfully Hurewicz, Menger, etc., go ``up" from $X_M$ to $X$.
On the other hand, it is not so clear what happens with a property like \textit{Alster}, since there can be expected to be compact subsets of $X$ that are not in $M$.

We may, instead of going from a countably closed elementary submodel up to $H_\theta$ or the entire universe,
go from the universe to an extension of it by countably closed forcing.
A typical argument is then that if some property involving the existence of a countable object holds in the extension,
it must have held in the original universe, since no new countable subsets of $V$ were added by the forcing.
Thus, if an open cover or sequence of open covers of $X$ in $V$ acquires some nice countable subcollection in the extension,
it must have had that nice subcollection already.
For example, if we find that $X$ is Lindel\"{o}f, Menger, Hurewicz, etc., in a countably closed forcing extension,
it must have had those properties to begin with.
Again, a property such as \textit{Alster} does not fit into this scheme,
because countably closed forcing does not in general preserve compactness, and moreover can adjoin new compact sets.

\bprb
If $X$ is $\sigma$-compact in a countably closed forcing extension, is it $\sigma$-compact?
\eprb

The analogous problem is also open for countably closed elementary submodels.

\bprb
If $M$ is a countably closed elementary submodel of $H_\theta$ for a sufficiently large regular $\theta$ with $X$ and its topology as members, then if $X_M$ is $\sigma$-compact, is $X$ also?
\eprb

For compactness, both problems have positive answers, and ``countably closed" is not needed.
This was noted earlier in the case of forcing; for elementary submodels, this was proved by L. R. Junqueira \cite{1J}.

%\nocite{*}
\bibliographystyle{amsplain}
\bibliography{NoteProdLind}

\newcommand{\noopsort}[1]{}
\providecommand{\bysame}{\leavevmode\hbox to3em{\hrulefill}\thinspace}
\providecommand{\MR}{\relax\ifhmode\unskip\space\fi MR }
% \MRhref is called by the amsart/book/proc definition of \MR.
\providecommand{\MRhref}[2]{%
  \href{http://www.ams.org/mathscinet-getitem?mr=#1}{#2}
}
\providecommand{\href}[2]{#2}
\begin{thebibliography}{10}

\bibitem{AAJT}
O.~T. Alas, L.~F. Aurichi, L.~R. Junqueira, and F.~D. Tall,
  \emph{Non-productively {L}indel\"of spaces and small cardinals}, Houston J.
  Math., in press.

\bibitem{Alster1988}
K.~Alster, \emph{On the class of all spaces of weight not greater than
  $\omega_1$ whose {C}artesian product with every {L}indel\"of space is
  {L}indel\"of}, Fund.\ Math. \textbf{129} (1988), 133--140.

\bibitem{1A}
\bysame, \emph{The product of a {L}indel\"of space with the space of
  irrationals under {M}artin's {A}xiom}, Proc.\ Amer.\ Math.\ Soc. \textbf{110}
  (1990), 543--547.

\bibitem{Arhangelskii1963}
A.~V. Arhangel'ski\u\i{}, \emph{On a class of spaces containing all metric
  spaces and all locally bicompact spaces}, Sov.\ Math.\ Dokl. \textbf{4}
  (1963), 751--754.

\bibitem{Arhangelskii1965}
\bysame, \emph{Bicompact sets and the topology of spaces}, Trans.\ Moscow
  Math.\ Soc. \textbf{13} (1965), 1--62.

\bibitem{Arhangelskii1986}
\bysame, \emph{Hurewicz spaces, analytic sets and fan tightness in function
  spaces}, Sov.\ Math.\ Dokl. \textbf{33} (1986), 396--399.

\bibitem{Arhangelskii2000}
\bysame, \emph{Projective $\sigma$-compactness, $\omega_1$-caliber, and
  {$C_p$}-spaces}, Topology Appl. \textbf{104} (2000), 13--26.

\bibitem{Aurichi}
L.~F. Aurichi, \emph{{$D$}-spaces, topological games and selection principles},
  Topology Proc. \textbf{36} (2010), 107--122.

\bibitem{AurichiTall}
L.~F. Aurichi and F.~D. Tall, \emph{{L}indel\"of spaces which are
  indestructible, productive, or ${D}$}, Topology Appl., to appear.

\bibitem{Barr2007}
M.~Barr, J.~F. Kennison, and R.~Raphael, \emph{Searching for absolute
  $\mathcal{CR}$-epic spaces}, Canad.\ J. Math. \textbf{59} (2007), 465--487.

\bibitem{1BT}
T.~Bartoszy\'nski and B.~Tsaban, \emph{Hereditary topological diagonalizations
  and the {M}enger-{H}urewicz conjectures}, Proc.\ Amer.\ Math.\ Soc.
  \textbf{134} (2006).

\bibitem{vanDouwen}
E.~K. {\noopsort{Douwen}}{van Douwen}, \emph{The integers and topology},
  Handbook of {S}et-theoretic {T}opology (K.~Kunen and J.~E. Vaughan, eds.),
  North-Holland, Amsterdam, 1984, pp.~111--167. \MR{MR776622 (87f:54008)}

\bibitem{vanDouwenPfeffer}
E.~K. {\noopsort{Douwen}}{van Douwen} and W.~F. Pfeffer, \emph{Some properties
  of the {S}orgenfrey line and related spaces}, Pacific J. Math. \textbf{81}
  (1979), 371--377.

\bibitem{Eisworth2007}
T.~Eisworth, \emph{On {$D$}-spaces}, Open {P}roblems in {T}opology {II}
  (E.~Pearl, ed.), Elsevier, Amsterdam, 2007, pp.~129--134.

\bibitem{Engelking1989}
R.~Engelking, \emph{General {T}opology}, Heldermann Verlag, Berlin, 1989.

\bibitem{GM84}
F.~Galvin and A.~W. Miller, \emph{$\gamma$-sets and other singular sets of real
  numbers}, Topology Appl. \textbf{17} (1984), 145--155.

\bibitem{1GN}
J.~Gerlits and Zs.\ Nagy, \emph{Some properties of ${C(X)}$, {I}}, Topology
  Appl. \textbf{14} (1982), 152--161.

\bibitem{Gruenhage2009}
G.~Gruenhage, \emph{A survey of {$D$}-spaces}, Set Theory and its Applications,
  Contemp.\ Math., {\rm ed.\ L. Babinkostova, A. Caicedo, S. Geschke, M.
  Scheepers, 2011}, pp.~13--28.

\bibitem{Hansell}
R.~W. Hansell, \emph{Descriptive topology}, Recent Progress in General Topology
  (M.~Hus\v{e}k and J.~van Mill, eds.), North-Holland, Amsterdam, 1992,
  pp.~275--315.

\bibitem{HenriksenIsbell}
M.~Henriksen and J.~R. Isbell, \emph{Some properties of compactifications},
  Duke Math.\ J. \textbf{129} (1957), 83--105.

\bibitem{HIJ}
M.~Henriksen, J.~R. Isbell, and D.~G. Johnson, \emph{Residue class fields of
  lattice-ordered algebras}, Fund.\ Math. \textbf{50} (1961), 107--117.

\bibitem{JayneRogers}
J.~E. Jayne and C.~A. Rogers, \emph{Borel isomorphisms at the first level --
  {I}}, Mathematika \textbf{26} (1979), 125--179, Borel isomorphisms at the
  first level: corrigenda et addenda. {\it Mathematika 27} (1980), 236--260.

\bibitem{1J}
L.~R. Junqueira, \emph{Upwards preservation by elementary submodels}, Top.\
  Proc. \textbf{36} (2010), 107--122.

\bibitem{Junqueira1998}
L.~R. Junqueira and F.~D. Tall, \emph{The topology of elementary submodels},
  Topology Appl. \textbf{82} (1998), 239--266.

\bibitem{1JW}
W.~Just and M.~Weese, \emph{Discovering {M}odern {S}et {T}heory, {II}}, Amer.\
  Math.\ Soc., Providence, 1997. \MR{1474727 (99b:03001)}

\bibitem{Michael1971}
E.~A. Michael, \emph{Paracompactness and the {L}indel\"of property in finite
  and countable {C}artesian products}, Compositio Math. \textbf{23} (1971),
  199--214.

\bibitem{Przymusinski1980}
T.~C. Przymusi\'nski, \emph{Normality and paracompactness in finite and
  countable {C}artesian products}, Fund.\ Math. \textbf{105} (1980), 87--104.

\bibitem{Przymusinski1984}
\bysame, \emph{Products of normal spaces}, Handbook of {S}et-{T}heoretic
  {T}opology (K.~Kunen and J.~E. Vaughan, eds.), North-Holland, Amsterdam,
  1984, pp.~781--826.

\bibitem{1R}
M.~Repicky, \emph{Another proof of {H}urewicz theorem}, Tatra.\ Mt.\ Math.\
  Publ., to appear.

\bibitem{1RZ}
D.~Repov\v{s} and L.~Zdomskyy, \emph{On the {M}enger covering property and
  {$D$} spaces}, Proc.\ Amer.\ Math.\ Soc., to appear.

\bibitem{Tall}
F.~D. Tall, \emph{{L}indel\"of spaces which are {M}enger, {H}urewicz, {A}lster,
  productive, or {$D$}}, Topology Appl., to appear.

\bibitem{1T}
B.~Tsaban, \emph{{M}enger's and {H}urewicz's {P}roblems: {S}olutions from
  ``{T}he {B}ook" and refinements}, Set Theory and its Applications, Contemp.\
  Math., {\rm ed.\ L. Babinkostova, A. Caicedo, S. Geschke, M. Scheepers, 2011,
  211--226}.

\bibitem{1TO}
B.~Tsaban and T.~Orenshtein, \emph{Linear $\sigma$-additivity and some
  applications}, Trans.\ Amer.\ Math.\ Soc. \textbf{363} (2011), 3621--3637.

\bibitem{SF1}
L.~Zdomskyy, \emph{A semifilter approach to selection principles}, Comm. Math.
  Univ. Carolinae \textbf{46} (2005), 525--539.

\end{thebibliography}

\end{document}